\documentclass[11pt]{article}

\usepackage{latexsym,amsmath,amssymb,theorem}
 \usepackage[applemac]{inputenc}  
\usepackage[all]{xy}
\usepackage{mathtools}

\usepackage{color} 
\usepackage{hyperref}

{\theorembodyfont{\itshape}
\newtheorem{Lemma}{Lemma}
\newtheorem{proposition}[Lemma]{Proposition}
\newtheorem{Fact}[Lemma]{Fact}

}
{\theorembodyfont{\upshape}

\newtheorem{Remark}[Lemma]{Remark}

}

\newcommand{\CT}{\mathcal{T}}
\newcommand{\CS}{\mathcal{S}}

\newcommand{\BB}{\mathcal{B}}
\newcommand{\BA}{\mathcal{A}}

\newcommand{\BS}{\mathcal{S}}
\newcommand{\BT}{\mathcal{T}}

\newcommand{\op}{\mathsf{op}}

\newcommand{\ob}{\mathsf{ob}}

\newcommand{\Eq}{\mathsf{Equ}}
\newcommand{\cod}{\mathsf{cod}}

\newcommand{\Set}{\mathit{Set}}

 \newcommand{\Equ }{\mathsf{Equ}}

\newcommand{\xra}{\xrightarrow}
\newcommand{\im }{\mathrm{Im}\,}

\begin{document}

\title{The  history of the General Adjoint Functor Theorem}
%

\author{Hans - E. Porst\\
 Besselstr. 65, D-28203 Bremen, Germany.\\  porst@uni-bremen.de}
 

\date{}
\maketitle

%

\begin{abstract} 
Not only motivated by the fact that the publication of the GAFT first appeared 60 years ago in print we reconstruct its history and so show that  it is no exaggeration to claim that it has appeared already 75 years ago!

\vskip 6pt
\noindent
 {\bf MSC 2020}:  {18-03, 18A40 }  \newline
\noindent {\bf Keywords:} {Universal morphism, left adjoint functor.  }

\end{abstract}

\section*{Introduction}

The General Adjoint Functor Theorem  appeared in print for the first time in Freyd's book \em Abelian Categories \em \cite{Freyd} in  1964; insiders had been aware of it earlier,  due to the fact that P. Freyd's PhD thesis (1960) containing it was apparently distributed  widely. Today  it is often called \em Freyd's Adjoint Functor Theorem\em, because it appeared under this name  in the first edition of Mac Lane's CWM \cite{MacL}. Consequently, one may call   2024 the year when this theorem turns 60.

In fact  --- and that is the topic of this note --- this theorem already turned 75 in 2023 and, moreover,  rather should be called the \em Samuel-Freyd Adjoint Functor Theorem\em. 
This suggestion may sound strange since  in 1948 the notion of an adjoint functor wasn't known at all.
However, in 1948 appeared P. Samuel's paper \em On universal mappings and free topological groups \em\cite{Samuel}. And this paper 
contains an existence theorem of what the author called \em universal mappings \em  --- and the core of the proof of the GAFT is the construction of universal arrows!

The author apparently was not yet familiar  with the concepts of  \em Category \em and \em Functor\em, introduced in the then recent 
Eilenberg-Mac Lane-paper \cite{EM2}. 
So he based his construction on an axiomatization of what later became known as \em  concrete categories \em  whose forgetful functors into  the category of sets create limits --- the first such axiomatization in the literature. 
It is somewhat irritating that in the later literature there seems to be no 
reference to this fact.

One probably should  recall in this context  that  the idea of the concept of a concrete category over the category of sets (though not that of concrete functors)  was in the air already in the 1930s 
(see Section \ref{tg}) and the annotation to Emmy Noether in Section \ref{sec:fin} below, as was the idea of universal maps (see Section \ref{tg}).


\section{The first examples of universal maps}

\subsection{Free groups (1920s)}

The first known example  of a universal map seems to be the embedding of the generating set $X$ into the free group $F(X)$ over $X$. This wasn't really seen as an important observation since it was  considered to be a property following quite obviously  from the construction of $F(X)$ as a group of (equivalence classes) of words in the letters $x$ and $x^{-1}$ with $x\in X$, \em ``free" \em from any relations but those of the group axioms. (The terminology \em free \em group was introduced by Nielsen \cite{Nielsen} in 1921). However, since the end of the 1920s algebraists started to use this characterization more often as it was considered more easy to use than the original construction.

\subsection{The Stone-$\check{\rm\bf{C}}$ech-compactification (1930/1937)}

 The first non-trivial example of a universal map seems to be 
what is called today the Stone-$\check{\rm {C}}$ech-compactification of a completely regular space. ${\rm\check{C}}$ech \cite{Cech} attributes the origins of this to Tychonoff \cite{Tychonoff} by  writing `` ... Tychonoff proves ... that, given a completely regular space $S$, there exists a bicompact Hausdorff space $\beta(S)$ such that (i) S is dense in $\beta(S)$, (ii) any bounded continuous real function defined in the domain S admits of a continuous extension to the domain $\beta(S)$". 
He then proves (certainly without using categorical language) that this result can be generalized to saying that the embedding of $S$ into  $\beta(S)$ is $E$-universal for the  functor $E$ embedding the category completely regular spaces into the category of compact Hausdorff spaces and that this property characterizes $\beta(S)$. {In contrast to the algebraists view on  the characterization  of free groups by their universal property, the topologists considered  this  as an  important result.}  
In ${\rm\check{C}}$ech's paper it also becomes clear that there was a size problem to be solved in order to get  a solution of the respective universal mapping problem.

\subsection{Free topological groups (1941 - 1945)}\label{tg}
As early as 1941 \cite{Markoff1} there appears the construction of a  universal map where at least its extended version \cite{Markoff2}, containing the motivation for this endeavour,  
 has a certain categorical flavour: Here Markov motivates his interest in free topological groups by the analogy between  the relations between 
 \em groups with group homomorphisms \em and \em sets with maps \em  on the one hand and \em topological groups with continuous homomorphisms \em and \em topological spaces with continuous maps \em  contrariwise. 
 He interprets free groups $F(X)$ over a set $X$ as given by the universal property of the embedding of the generating set $X$ into $F(X)$ and, hence, asks for an analog property of the embedding of a ``generating" \em topological space \em  $X$ into a  \em topological group \em $F(X)$.
Since only Hausdorff topological groups were seen to be of interest (and these are completely regular automatically), and he wanted the generating space $X$ to be a subspace of the respective topological group $F(X)$ in analogy to the groups case, he assumed $X$ to be completely regular. 
He introduces sets of \em multinorms \em in order to solve the size problem. 

Motivated by \cite{Markoff1} Kakutani \cite{Kakutani} gave a much simpler proof of the existence of free topological groups by solving the size problem in showing that, for any completely regular space $X$, the set of all topological groups $G$ with $card (G)\leq max(card(X),\aleph_0)$ is a solution set and then using the by now standard construction of a universal embedding of $X$ into $F(X)$ only a couple of years later.

\section{Samuel's universal mapping problem (1948)}\label{sec:2}

For  Samuel, who already had successfully used an analog of $\check{\rm {C}}$ech's solution of the size problem in his PhD thesis  \cite{SamuelTh}, the appearance of Markov's paper \cite{Markoff2} seems to have been the decisive  motivation to look for a more general approach to the constructions mentioned above and made him  write  the paper \cite{Samuel}. 

\subsection{Concrete categories}\label{sec:2.1}

Not being aware of the then recent paper \cite{MacL-ori} 
he did not use the notions \em category \em and \em functor\em, but rather abstracts his axioms from Bourbaki's notion of  \em structure\em.  
He, hence, was talking about $T$-sets and $T$-mappings for some ``kind of structure" $T$.  
This is in some detail:  
\begin{enumerate}
\item Given a ``kind of structure" $T$ then, for each set $Y$ there exists a set $T(Y)$ of ``$T$-structures on $Y$".
\item  $T$-sets are pairs $(Y,\tau_Y)$ where $Y$ is a set and $\tau_Y \in T(Y)$.
\item $T$-mappings  $(Y,\tau_Y)\xra{} (Y',\tau_{Y'})$ are (certain) maps $Y\xra{}Y'$.
\item $T$-isomorphisms  $(Y,\tau_Y)\xra{} (Y',\tau_{Y'})$ are (certain) bijective $T$- mappings. 
\item For each $T$-set $(Y,\tau_Y)$ the identity $id_Y$ is a $T$-isomorphism.
\end{enumerate}
He then assumes that $T$-sets and $T$-mappings satisfy certain axioms. 
It is obvious, however, of how to translate this  into the language of category theory  --- more precisely, into the language of concrete categories over the category $\Set$ of sets in the sense of \cite{AHS}. 
The original formulation of his axioms can be found in the Appendix, where also some remarks are added concerning this translation.

\begin{description}
\item[The category $\CT$  of $T$-sets:] 
Given a kind of structure $T$, the {\it $T$-sets} $(Y,\tau)$ 
and {\it $T$-mappings}  $(Y,\tau)\xra{f} (Y',\tau')$ 
%
%
%
 form a concrete category $\CT$ whose underlying functor $\CT\xra{|-|}\Set$ 
acts  as  $(Y,\tau)\mapsto Y$. 
 \end{description}

\begin{description}
\item[Properties of $\CT$:] 
The forgetful functor $\CT\xra{ |-|}\Set$ creates products, equalizers, and intersections and reflects isomorphisms. 
 \end{description}
 One, hence,  can form  for each $T$-object $(Y,\tau)$ and any subset $Z\subset Y$ the intersection of  all  subobjects of $(Y,\tau)$ containing $Z$. This is the smallest  subobject $(\bar{Z},\bar{\zeta})$  of $(Y,\tau)$ with $Z\subset\bar{Z}$,  called the \em $T$-closure \em  of  $Z$.
\begin{description} 
\item[Axiom S$_3$] 
There is a monotonic\footnote{Monotonicity  of $\kappa_T$ isn't explicitly mentioned by Samuel, but tacitly assumed. } function of cardinals $\kappa_T$ such that, for any  $T$-set $(Y,\tau)$ and  for any subset $Z\subset Y$, $card \bar{Z} \leq  \kappa_T(card Z)$.
\end{description}
In particular the following holds:
\begin{Fact}\label{fact:}\rm
For every  $\CT$-object $(Y,\tau)$ and any subset $ Z\subset Y$ one has the inequalities $card \bar{Z}\leq card Y =card |(Y,\tau)|\leq \kappa_T (card Y)$ such that, in particular, the (isomorphism classes of) all $\CT$-objects over a set $Y$ form a set.
\end{Fact}

A paradigmatic example of a structure $T$ satisfying  Axiom S$_3$ is the structure $T_{Grp}$ of groups:  If $G$ is a group then, for any subset $X$ of $G$, its closure $\bar{X}$ is the subgroup $\langle X\rangle$ generated by $X$ and $card \langle X\rangle$  clearly is smaller than the cardinality of $\bigcup_n \{ (x_1,\ldots x_n)  \mid x_i\in X\uplus X\}  =: \kappa_{T_{Grp}}(X)$.

\begin{description}
\item[{\em S-T\em-mappings} \rm{(Axioms S-T)}: ] Given structures $S$ and $T$ then, for each  object $(X,\sigma)$ in $\CS$ and   each object $(Y,\tau)$ in $\CT$ there is given a set
  $M((X,\sigma),(Y,\tau))$ of  {\it S-T-mappings}, i.e., of maps $X\xra{\phi}Y$   subject to the following condition:
\begin{description}
\item[\rm{(S-T)$_1$}]  \   If  $(X,\sigma)\xra{\phi'}(Y',\tau')$ is an  \em S-T\em-mapping and $(Y',\tau')\xra{f}(Y,\tau)$ is a morphism in $\CT$, 
 then $\phi = (X,\sigma)\xra{\phi'}(Y',\tau')\xra{f}(Y,\tau)$  is an  \em S-T\em-mapping.
\item[\rm{(S-T)$_2$}] \  For any family of \em S-T\em-mappings $(X,\sigma)\xra{\phi_i}(Y_i,\tau_i)$ the map  $(X,\sigma)\xra{\phi}\Pi_i(Y_i,\tau_i)$ induced by this family by the product property (in $\Set$) is an \em S-T\em-mapping. 
\end{description}
 \end{description}

%

\begin{Remark}\label{rem:miss}\rm
There is, however, an axiom missing. In the first step of his proof of the solution of the universal mapping problem Samuel uses the following fact, very much in line with his Axiom I$_2$:
\end{Remark}

\begin{Fact}\rm 
If $(X,\sigma)\xra{\phi}(Y,\tau)$ is an \em S-T\em-mapping and $\overline{\im\phi}$ the closure of the image of $\phi$, i.e.,   the intersection of all subobjects $(Y',\tau')$ of $(Y,\tau)$ such that $\phi$ factors (in $\Set$) as $X\xra{\phi'}Y'\hookrightarrow Y$,  then the contraction (corestriction) $(X,\sigma)\xra{\phi'}\overline{\im\phi}$ of $\phi$ is an  \em S-T\em-mapping. 
\end{Fact}
This obviously requires the following additional axiom :
\begin{quote}
\rm{ (S-T$)_3$:}   If $(X,\sigma)\xra{\phi}(Y,\tau)$ is  an  \em S-T\em-mapping and 
$(Y',\tau')\xra{i}(Y,\tau)$ is the embedding of a $T$-subobject such that $\phi$ factors  (in $\Set$) as $\phi = X\xra{\phi'}Y'\xra{i}Y$
then $(X,\sigma)\xra{\phi'}(Y',\tau')$ is an  \em S-T\em-mapping. 
\end{quote}
This condition is satisfied an all examples given in \cite{Samuel}.  
Note that in view of all other axioms this would be equivalent to 
\begin{quote} \rm{(S-T)$_3'$:} \ 
 Let  $(Y_i,\tau_i)\xra{m_i}(Y,\tau)$ be a family  of subobject embeddings $(Y_i,\tau_i)$ and $(\bar{Y},\bar{\tau}):=\bigcap_i(Y_i,\tau_i)\xra{m}(Y,\tau)$ the embedding of its intersection. Then, for any \em S-T\em-mapping $(X,\sigma)\xra{\phi}(Y,\tau)$ and any family of mappings $X\xra{\phi_i}Y_i$ with $X\xra{\phi_i }Y_i\xra{m_i} Y =\phi $,  the map    $X\xra{\psi}\bar{Y}$ with $X\xra{\psi}\bar{Y}\xra{m}Y =\phi$ resulting from the intersection property  is an  \em S-T\em-mapping  $(X,\sigma)\xra{\psi}(\bar{Y},\bar{\tau})$.
\end{quote}

\begin{Remark}\label{rem:2}\rm
 The concept of \em S-T\em-mappings can be seen as an attempt to express a relation between two such categories $\CS$ and $\CT$ in the absence of the concept of functor. Consequently, using 
categorical language, this is  best and most naturally   done by thinking of a concrete functor $E\colon \CT\xra{ }\CS$, that is, a functor from $\CT$ to $\CS$ making the following diagram commute
\begin{equation*}\label{diag:}
\begin{aligned}
\xymatrix@=2em{\CT
{ }\ar[rr]^{ E} \ar[dr]_{|-|}& &\CS \ar[dl]^{|-| }     \\
& {\Set }     
}
\end{aligned}
\end{equation*}
 and choosing   the
 sets $M((X,\sigma),(Y,\tau))$   to be  the hom-sets $\CS((X,\sigma),E(Y,\tau))$\footnote{For alternative interpretations see Section \ref{sec:fin}.} 
This interpretation is supported by Axiom (S-T)$_1$ and applies to all examples mentioned in \cite{Samuel}. Axioms (S-T)$_2$ and (S-T)$_3$ then are satisfied iff the  functor $E$ preserves limits (hence products and intersections in particular) automatically, since the respective forgetful functors are assumed to create limits.
\end{Remark}

\subsection{The problem of universal mappings}

 When translated into categorical  language as suggested in Remark \ref{rem:2}, the \em problem of universal mappings \em  then is the following.  
 \begin{quote}
 Let $\CT\xra{E}\CS$ be a concrete functor.
    Given an object $(X,\sigma)$ in $\CS$, find an   object $(Y_0, \tau_0)$ in $\CT$  and an  $\CS$-morphism $\psi\in  \CS((X,\sigma),E(Y_0,\tau_0))$ such that for every  $\CS$-morphism $\phi\in  \CS((X,\sigma),E(Y,\tau))$ there exists a unique $\CT$-morphism $(Y_0, \tau_0)\xra{f}(Y,\tau)$ with $\phi= E(f)\circ\psi$.
\end{quote}

Samuel's crucial (though simple) observation was that this problem essentially is a size problem: if  the category $\CT$ were small one could, since  $E$ preserves products, simply  choose, with  $\Delta_{(X,\sigma)} = \{ (X,\sigma)\xra{\phi}E(Y,\tau)  \mid  (Y,\tau)\in \ob\CT\}$,\vspace{-1ex}
$$(Y_0,\tau_0) :=\prod_{\ob\CT}(Y,\tau)^{\CS((X,\sigma),E(Y,\tau))} 
\text{\ \ such that \ \ }E(Y_0,\tau_0) =\prod_{\Delta_{(X,\sigma)}}\cod\phi\vspace{-2ex}$$
and 
$\CS(X,\sigma)\xra{\psi} E(Y_0,\tau_0)$
 to be the unique $\CS$-morphism induced by the set of all $\CS$-morphisms $\phi\in \Delta_{(X,\sigma)}$ 
 by the product property, that is, the unique map with $E(\pi_\phi)\circ \psi =\phi$.
Since usually $\CT$ fails to be small one is lead to the question whether one can replace in the definition of $(Y_0,\tau_0)$ above the index set $\ob\CT$ by a subset of the class $\ob\CT$.

Samuel's  idea to overcome this problem explains why he introduced 
 Axiom S$_3$\footnote{ See the last paragraph of this section for  a possibly more suitable formulation of this axiom.}: This allows him to 
 replace,  for any $\CS$-object $(X,\sigma)$, $\ob\CT$ by a representative set $\Lambda_{(X,\sigma)} $ of the class of all $(Z,\zeta)\in \ob\CT $ with $    card(Z)\leq \kappa_T(card X) $   and, accordingly, the set $\Delta_{(X,\sigma)}$ by 
 $\{\lambda\in \CS((X,\sigma),E(Z_\lambda,\zeta_\lambda ))\mid (Z_\lambda,\zeta_\lambda )\in \Lambda_{(X,\sigma)} \}$.
  
The first step of his construction then is the following:
If $(X,\sigma)\xra{\phi}E(Y,\tau)$ is an $\BS$-morphism, and  
$\phi= (X,\sigma)\xra{\phi'}E(\bar{Y'},\upsilon) \xra{ E(e_\phi) }E(Y,\tau) $ with $Y' =\overline{\im\phi}$ and $e_\phi$ the embedding of its closure,   then   $\phi'$ is an $\CS$-morphism  by Axiom (S-T)$_3'$ (respectively, the fact that $E$ preserves intersections) such that  (by monotonicity of $\kappa_T$) $card(\bar{Y}')\leq   \kappa_T(card(\im\phi) )\leq \kappa_T(card X)$  and, hence,  (up to an isomorphism) 
   $ (\bar{Y'},\upsilon) = (Z_{\lambda_\phi},\zeta_{\lambda_\phi}) \in \Lambda_{(X,\sigma)}$.
Forming the product $(Y'_0,\tau'_0):= {\prod_{\Lambda_{(X,\sigma)}}} (Z_{\lambda_\phi},\zeta_{\lambda_\phi}) $ one so obtains  by Axiom \rm{(S-T)$_2$} (respectively,  the fact that $E$ preserves products) an  $\CS$-morphism 
  $(X,\sigma)\xra{\psi'} E(Y'_0,\tau'_0)=  {\prod_{\Lambda_{(X,\sigma)}}} E(Z_{\lambda_\phi},\zeta_{\lambda_\phi}) $ such that the following diagram commutes, where $\pi_{i_\phi}$ is the respective projection.
\begin{equation*}\label{diag:} 
\begin{aligned}
\xymatrix@=2.5em{
{ (X,\sigma) }\ar[r]^{\psi' } \ar[d]_\phi\ar[dr]_{\phi'}& E(Y'_0,\tau'_0) \ar[d]^{ E(\pi_{i_\phi)}}     \\
E(Y,\tau)&  E(Z_{\lambda_\phi},\zeta_{\lambda_\phi})  \ar[l]^{E(e_\phi)}     
}
\end{aligned}
\end{equation*}

Thus, what remains to be done is to look for a subobject  embedding $i\colon (Y_0, \tau_0)\xra{}(Y'_0,\tau'_0)$ such that for any pair of $\CT$-morphisms $f', g'\colon(Y'_0,\tau'_0)\xra{} (Y,\tau)$ with $E(f')\circ\psi' = E(g')\circ\psi'  = \phi$ the $T$-set 
 $(Y_0, \tau_0)$  is a subobject of its  equalizer $\Eq(f',g')$ with embedding $i'$.   Then, by   the assumption that $E$  preserves equalizers  
$\psi'$ would  factor, with $j'$ the equalizer embedding, as 
\begin{equation*}\label{eqn:}
 \begin{aligned}
\xymatrix@=3em{
{ \psi' = E(X,\sigma)}\ar[r]^{\psi } 
& E(Y_0, \tau_0)\ar@/_1.5pc/@{->}[rr]_{E(i)}\ar[r]^{E(i') }& E(\Eq(f',g'))\ar[r]^{E(j')} & E(Y'_0, \tau'_0)
}
\end{aligned}
\end{equation*}
for any pair of $\CT$-morphisms $f',g'\colon (Y'_0, \tau'_0)\xra{}(Y,\tau)$ satisfying the equations  $\phi = E(f')\circ\psi' = E(g')\circ\psi'$. 
Consequently the $\CT$-morphism $f:= f'\circ i   = g'\circ i $ would be the only $\CT$-morphism satisfying the equation $E(f) \circ\psi = \phi$, such that the $\CT$-object $(Y_0, \tau_0)$ together with the $\CS$-morphism $\psi$ would be a solution of the universal mapping problem.

This is what Samuel does in the second step as follows: He chooses $(Y_0, \tau_0)$ to be the $T$-closure of $\im\psi'$.
 Since $\im\psi'$ is contained in the equalizer of every pair $(E(f'),E(g'))$ satisfying 
$E(f')\circ\psi' = E(g')\circ\psi'  = \phi$ for some $\phi$,  and this equalizer is $E(\Equ(f',g' ))$ where $\Equ(f',g') $ is the equalizer of $(f',g')$ in $\CT$ (since $\CT$ has equalizers and these are preserved by $E$), $(Y_0, \tau_0)$ is contained in all of these equalizers as well. By the above the morphism   $ (X,\sigma)\xra{\psi}E (Y_0, \tau_0) $ is the desired solution. 
The situation is illustrated by the following commutative diagram, where  $f:= e_{\phi'}\circ \pi_{{\phi'}}\circ i$ now is the unique $\CT$-morphism $(Y_0,\tau_0)\xra{} (Y,\tau)$ with $\phi = E(f)\circ\psi$.

\begin{equation*}\label{diag:}
\begin{aligned}
\xymatrix@=2em{
{(X,\sigma)}\ar[r]^{\!\!\!\!\psi }\ar[d]^{\phi} \ar[drr]_{\phi'}\ar@/^1.5pc/@{->}[rr]^{\psi'}
& 
E(Y_0, \tau_0)\ar[r]^{\!\!\!\!\!\!\!\!\!\!\!\!\!\!\!\!\!\!\!\!\!\!\!\!\!\!\! E(i)}& E(Y'_0,\tau'_0)
=\smashoperator{\prod_{\Lambda_{(X,\sigma)}}} E(Y_\phi,\tau_\phi)
 \ar[d]^{{E(\pi_{\phi')}} }    \\
E(Y,\tau) &  &\ar[ll]^{\!\!\!\!\!\!\!\!\!\!\!\!\!\!\!\!\!\!\! \ E(e_{\phi'})}  E(Y_{\phi'},\tau_{\phi'})   
}
\end{aligned}
\end{equation*}

In categorical language  the above result is the following

\begin{proposition}\label{prop:1}
Let $\BS$ and $\BT$ be concrete categories over $\Set$ whose forgetful functors creates limits. Then, given  a   concrete and, hence, limit preserving functor $E\colon  \CT\xra{}\CS$,   for each $\BS$-object $(S,\sigma)$ there exists an $E$-universal arrow $(S,\sigma)\xra{ } E(T,\tau)$ 
provided that $\CT$ satisfies Axiom S$_3$. 
\end{proposition}

Looking at the proof above one sees that Axiom S$_3$  is only used to be able to define the set $  \Lambda_{(X,\sigma)}$ 
and make the following true:
\begin{itemize}
\item for each $\CS$-morphism $ (X,\sigma)\xra{\phi} E(Y,\tau)$ there exists some $\CS$-morphism $ (X,\sigma)\xra{\psi} E(Z,\zeta)$  and some $\CT$-morphism $(Z,\zeta)\xra{f}(Y,\tau) $ with  $(Z,\zeta)\in \Lambda_{(X,\sigma)}$ satisfying the equation $E(f)\circ\psi = \phi$.
\end{itemize}
In other words, it would have been sufficient (and more appropriate) 
instead of adding Axiom S$_3$ to the list of axioms of the categories of $T$-sets, to rather add the following axiom to the list of axioms of \em S-T\em-mappings:
\begin{description}
\item[SolSet ] For each $S$-set $(X,\sigma)$ there exists a set   $  \Lambda_{(X,\sigma)}$ of $T$-objects such that for each \em S-T\em-mapping $ (X,\sigma)\xra{\phi} (Y,\tau)$ there exists some \em S-T\em-mapping $(X,\sigma)\xra{\psi} (Z,\zeta)$  and  and some $T$-mapping  $(Z,\zeta)\xra{f}(Y,\tau) $ with  $(Z,\zeta)\in \Lambda_{(X,\sigma)}$  such that the following diagram commutes
\begin{equation*}\label{diag:}
\begin{aligned}
\xymatrix@=2em{
{ (X,\sigma) }\ar@{.>}[r]^{\psi } \ar[d]_{\phi}& (Z,\zeta) \ar@{.>}[dl]^{f }     \\
 { (Y,\tau) }     &
}
\end{aligned}
\end{equation*}
\end{description}


\section{Freyd's GAFT (1964)}

Proposition \ref{prop:1} immediately leads to the question of how to generalize it to more general limit preserving functors $G\colon\BA\xra{}\BB$.  By inspection one sees that what  is used in Samuel's proof  are the following statements. 

\begin{enumerate}
\item  The category $\BT$ has products, equalizers (and intersections) and these are preserved by $E$. 
(The assumption that $\CT$ and $\BS$ are concrete categories with  forgetful functors which lift limits is only used to get this.) 
\item  Condition {\bf SolSet} is satisfied.
\end{enumerate}

Obviously condition {\bf SolSet} implies Freyd's solution set condition for the functor $E$ and is used only  to get this.
Thus, by simply translating Samuel's solution of the universal mapping problem into the language of categories, functors and limits as done above
one essentially obtains Freyd's  GAFT \cite{Freyd};  (the trivial conclusion, that the existence of a $G$-universal map for each $\BB$-object $B$  leads to the existence of a left adjoint of $G$  is already contained in \cite{Kan}).
This --- together with the fact that Samuel was the first to not only look for \em analogies \em between the various  constructions of universal maps known in the early 1940s but for a \em common pattern \em of those --- motivates the suggestion of the introduction to rather call this theorem the \em Samuel-Freyd Adjoint Functor Theorem\em.

Without much doubt one may say that Kan already would have included this theorem in  \cite{Kan}, if he had been aware of Samuel's paper.

\section{Final remarks}\label{sec:fin}
\subsection{References and later interpretations}\label{inter}
Mac Lane knew  Samuel's paper already from reviewing it (see \cite{MacL-MR}). He here already observed that axiom \rm{ (S-T)$_3$} is missing in Samuel's  list of axioms; somewhat surprisingly, however, he does not refer to \cite{EM2}. 
He later refers to Samuel's  paper quite often: In the notes to \cite[\em Notes \em to Chapter III.]{MacL}  he  writes \em ``The bold step of really formulating the general notion of a universal arrow was taken by Samuel in 1948"; \em   in \cite[p. 341]{MacL-ori2} he even writes \em  ``In [1948], Samuel described universal constructions". \em
In \cite[p. 56]{MacL-ori} he gives an interpretation of Samuel's paper different from ours in that he assumes that the \em S-T \em sets  $M((X,\sigma),(Y,\tau))$ define a profunctor from $\CT$ to $\CS$, that is, a functor $\CS^\op\times\CT\xra{}\Set$. 

Ellermann \cite{El} interprets  Samuel's approach in a way similar to Mac Lane. Referring to Mac Lane's references to Bourbaki\footnote{These are  indirect references to Samuel (see  \cite{El}).} in the notes to \cite[\em Notes \em to Chapters IV. and V.]{MacL} he states inter alia 
\em``Samuel may have been even closer \em (to adjunctions) \em than Mac Lane surmised"\em, which is in view of the above more than true!

\subsection{Some historical remarks}\label{hist}

\begin{description}
\item[Markov,]more precisely,  A. A. Markov, Jr., worked in Saint Petersburg (then called Leningrad), which was for  most of the time between appearance of the papers \cite{Markoff1} and \cite{Markoff2} under siege  during WW II. This may explain the unusual fact that Markov published two papers with identical title (first the announcement of the results and later the proofs) with a four years time interval.
\item[Samuel]  continued his studies, after an interruption  due to WW II,  in Princeton where he got his PhD in 1947; before returning to France he submitted the paper \cite{Samuel}. In the early 1950s he became a member of the Bourbaki group and so his solution of the universal mapping problem was incorporated into \cite{Bourbaki};  unfortunately, Bourbaki then already had decided not to use the language of categories and the categorical flavour of \cite{Samuel} got lost.
\item[Mac Lane]  became  influenced by Emmy Noether during his two years stay in Göttingen in the early 1930s. He later occasionally called her ``grandmother of category theory", probably   since she was the first to routinely talk about \em groups \em and their \em homomorphisms\em,  \em rings \em and their \em homomorphisms\em, or  \em modules \em and their \em homomorphisms\em, in other words,  the (concrete) categories of the respective structures (see e.g. \cite{Noether} --- before (with few exceptions) only isomorphisms had been considered.
\end{description}

\section*{Appendix}
Samuel starts writing
\begin{quote}
{\small
Given a set E it is possible to
define on it certain kinds of structures, that is structure of ring, field, a
topological space. We shall denote by $S$ or $T$ certain kinds of structures. A set with a structure $T$ will be called a \em T-set\em: if $T$ is the structure of group the $T$-sets are the groups. }
\end{quote}
The axioms in \cite{Samuel} then are formulated as follows: 
\vspace{0.5em}

{\small
{\bf $T$-mappings:}  Given a kind of structure $T$ it happens very often that, for every pair $(E_1,E_2)$ of $T$-sets, there has been defined a family of mappings of $E_1$ into $E_2$ satisfying the following axioms:
\begin{description}
\item[A$_1$.\ ] \   Every $T$-isomorphism is a $T$-mapping.  
\item[A$_2$.\ ] \ If $f\colon E_1\xra{ }E_2$ and $g\colon E_2\xra{ }E_3$ are $T$-mappings, then the composite
mapping $g\circ f\colon E_1\xra{ }E_3$ is a $T$-mapping.
\item[A$_3$.\ ] A necessary and sufficient condition for a one-to-one mapping $f$ of $E_1$ onto $E_2$ to be a $T$-isomorphism is that $f$ and $f^{-1}$ be $T$-mappings. 
\end{description}

EXAMPLE. If $T$ is the structure of group the $T$-mappings are the homomorphisms; if $T$ is the structure of topological space the 
$T$-mappings are the continuous ones.

\vspace{0.5em}
{\bf Induced structures:}
Let now $\sigma$ and $\sigma'$ be two structures $T$ defined on $E$ and $E'\subset E$ 
respectively. We shall say that $\sigma'$ is induced by $\sigma$ when: 
\begin{description}
\item[I$_1$.\ ] \ The injection of $E'$ into $E$ is a $T$-mapping.
\item[I$_2$.\ ]  If $f\colon F\xra{ }E$ is a $T$-mapping and if $f(F)\subset E'$, then $f$ considered $f$
as mapping of $F$ into $E$ is a $T$-mapping.
\end{description}
If $E'\subset E$  is capable of an induced structure we shall say that $E'$ is $T$-closed. We suppose that the following axioms hold :
\begin{description}
\item[S$_1$.\ ] \  A subset of $E$ composed of all the elements where a family of
$T$-mappings takes the same value is $T$-closed.
\item[S$_2$.\ ]   Any intersection of $T$-closed sets is $T$-closed.
\item[S$_3$.\  ]  Cardinal $(\bar{E'}) \leq $ certain function of cardinal $(E')$, a function which depends only on the structure $T$.
\end{description}

{\bf Axioms for the cartesian products.} In many important cases it is possible, given a family $(E_\alpha)$ of $T$-sets, to define on the cartesian product $\prod_\alpha E_\alpha$, a structure $T$ which satisfies the following conditions:
\begin{description}
\item[P$_1$.\ ]    The projections (on the components) are $T$-mappings.
\item[P$_2$.\ ]   If the $f_\alpha\colon X\xra{ }E_\alpha$ are $T$-mappings, the product mapping $f\colon E \xra{ }\prod_\alpha E_\alpha$  (defined by $f(x) = (f_\alpha(x))$ is a $T$-mapping.
\end{description}

Given two kinds of structures $S$ and $T$, suppose we have defined the $T$-mappings, and also mappings of $S$-sets into $T$-sets, called the (\em S-T\em)-mappings, denoted by greek letters, and satisfying:
\begin{description}
\item[(S-T)$_1$.]The composite mapping $f \circ\phi$ of an (\em S-T\em)-mapping $\phi$ and of
a $T$-mapping $f$ is an \linebreak 
\mbox{\ \ \ } (\em S-T\em)-mapping.
\item[(S-T)$_2$.]The product mapping of a family of (\em S-T\em)-mappings is an
(\em S-T\em)-mapping.
\end{description}
}

%

\begin{Remark}\label{rem:4}\rm  
These axioms imply that  $T$-sets and $T$-mappings form a complete concrete category  $\CT$, with  limits created by the forgetful functor.
Indeed, the forgetful functor creates products  by axioms {\bf P}$_1$ and {\bf P}$_2$ and it creates equalizers and intersections by Axioms {\bf S}$_1$ and {\bf S}$_2$   in connection with Axiom {\bf I}$_2$.
Hence, isomorphisms  are reflected by the forgetful functor  $|-|$ (see e.g \cite[17.14.]{AHS}) 
as stated in Section \ref{sec:2.1}.  Conversely,  all such concrete categories satisfy the above axioms (with Axioms {\bf I} restricted to embeddings of equalizers and intersections). 

One, hence, well might say --- though he wasn't defining categories --- Samuel was the first to axiomatize a certain class of categories, namely those concrete categories over $\Set$ whose underlying functor creates limits, 
and this for a good reason.
\end{Remark}

\pagebreak


\end{document}